\documentclass{siamltex}
\usepackage{amsmath, amsfonts,
  verbatim, fancyhdr, graphicx, array, subfigure,
  color, listings, afterpage, epsfig, boxedminipage, rotating, url,
  grffile}

\newcommand{\citet}{\cite}
\newcommand{\citep}{\cite}

\newcommand{\matcol}{\ensuremath{\text{:}}}
\newcommand{\matcom}{\ensuremath{\text{,}}}

\newcommand{\nfro}[1]{\ensuremath{\|#1\|_{\text{F}}}}
\newcommand{\nmax}[1]{\ensuremath{\|#1\|_{\text{max}}}}
\newcommand{\ntwo}[1]{\ensuremath{\|#1\|_2}}
\newcommand{\nrmone}[1]{\ensuremath{\|#1\|_1}}
\newcommand{\nnz}{\ensuremath{\text{nnz}}}
\newcommand{\hmatrix}{{\it H}-matrix}

\title{H-Matrix and Block Error Tolerances}

\author{Andrew M.~Bradley\thanks{Dept.~of Geophysics, Stanford University ({\tt
      ambrad@cs.stanford.edu}). This research was supported by NSF grant
    EAR-0838267 and USGS grant G10AP00009-001.}}

\begin{document}
\sloppy
\maketitle

\begin{abstract}
We describe a new method to map the requested error tolerance on an \hmatrix{}
approximation to the block error tolerances. Numerical experiments show that the
method produces more efficient approximations than the standard method for
kernels having singularity order greater than one, often by factors of $1.5$ to
$5$ and at a lower computational cost.
\end{abstract}



\pagestyle{myheadings}
\thispagestyle{plain}
\markboth{A.M.~BRADLEY}{H-MATRIX AND BLOCK ERROR TOLERANCES}

\section{Introduction}
A completely dense matrix $B$ arising from an integral or sum involving a
singular kernel over a surface or particles can be approximated efficiently by a
\emph{hierarchical} matrix, called an \emph{\hmatrix{}} \citep{hlib}. Let $\bar
B$ be the \hmatrix{} approximation to $B$.

Let the surface be discretized by $N$ elements or let there be $N$ particles.
The two cases differ essentially only in how a matrix entry of $B$ is
calculated. If $B$ includes all and only pair-wise interactions, then $B \in
\mathbb{R}^{N \times N}$. In this paper $B$ is always square, but the results
extend to the rectangle case $B \in \mathbb{R}^{M \times N}$ by replacing
occurrences of $N^2$ by $MN$. The procedure to construct an \hmatrix{} has four
parts. First, a cluster tree over the particles is formed. The cluster tree
induces a (nonunique) symmetric permutation of $B$. For notational brevity,
hereafter we assume $B$ is already ordered such that the identity matrix is a
permutation matrix induced by the cluster tree. Second, pairs of clusters are
found that satisfy certain criteria; associated with such a pair is a
\emph{block} of $B$. Third, the requested error tolerance $\varepsilon$
(hereafter usually just \emph{tolerance}) is mapped to block tolerances. A
tolerance specifies the maximum error allowed. Fourth, each block is
approximated by a \emph{low-rank approximation} (LRA) that satisfies that
block's tolerance. A number of LRA algorithms are available.

An LRA to a block $B_i$ can be efficiently expressed as an outer product of two
matrices $U$ and $V$: $B_i \approx \bar B_i = U V^T$. Let $r$ be the number of
columns in $U$; then $r$ is the maximum rank of $\bar B_i$ and the rank if $U$
and $V$ have independent columns, as is always the case in this context. Let
$\nnz(B_i)$ be the number of nonzero numbers required to represent $B_i$
exactly, and similarly for other matrices. Hence $\nnz(B) = N^2$; if $B_i \in
\mathbb{R}^{m \times n}$, then $\nnz(B_i) = m n$ and $\nnz(\bar B_i) = \nnz(U) +
\nnz(V) = (m + n) r$; and if $\bar B$ is partitioned according to the index set
${\cal I}$, then $\nnz(\bar B) = \sum_{i \in {\cal I}} \nnz(\bar B_i)$. The
efficiency of an \hmatrix{} approximation $\bar B$ is measured by its
\emph{compression}, $\nnz(B) / \nnz(\bar B) = N^2 / \nnz(\bar B)$. A standard
application of an \hmatrix{} approximation is to a sequence of matrix-vector
products (MVP) with dense vectors. The work to compute an MVP is proportional to
$\nnz(\bar B)$.

In this paper we focus on the third part of constructing an \hmatrix{}: how to
map the requested tolerance $\varepsilon$ on $\bar B$ to block tolerances. An
efficient method has at least these three properties. First, the tolerance is
met. Second, the tolerance is not exceeded by much: any unrequested additional
accuracy reduces the compression. Third, the method must not add much work to
the construction process relative to a standard method; indeed, if the resulting
compression is greater, the work will likely be less. We propose a method based
on the inequality \eqref{eq:mwrem} that has these properties.

It may be worth emphasizing the importance of an algorithm's meeting a tolerance
without exceeding it by much. Many computational mathematical algorithms trade
between speed and accuracy. A tolerance has two roles: first, it is the means by
which the user controls the tradeoff; second, the associated error bound, either
by itself or in combination with those from other algorithms in a larger
framework, gives the accuracy of the approximate solution to the problem. A user
who detects that an algorithm is giving much greater accuracy than requested is
tempted to request less accuracy. Such a procedure requires more work in the
form of repeated solutions of trial problems until the speed-accuracy tradeoff
\emph{appears}---but cannot be known with certainty---to be satisfactory; and it
sacrifices an assured error bound because the relationship between tolerance and
error bound is not clear. The method we shall describe does not principally
increase compression, relative to the standard method, for a given (unknown)
\emph{achieved} error (though it does in important cases in our experiments);
rather, it increases compression for a given \emph{requested} error tolerance by
achieving an error that is not much less than the tolerance.

\section{The meaning of the requested error tolerance}
In this paper and as often implemented in \hmatrix{} software, the tolerance
$\varepsilon$ specifies the maximum allowed Frobenius-norm-wise relative error
of the approximation $\bar B$ to $B$. Let $E \equiv B - \bar B$. It is requested
that $\bar B$ satisfy
\begin{equation} \label{eq:fnwre}
  \nfro{E} \le \varepsilon \nfro{B} \quad \text{or, equivalently,} \quad
  \frac{\nfro{B - \bar B}}{\nfro{B}} \le \varepsilon.
\end{equation}
If \eqref{eq:fnwre} holds, then $\nfro{E x} \le \nfro{E} \ntwo{x} \le
\varepsilon \nfro{B} \ntwo{x}$. Rearranging, the relative error in an MVP is
bounded as
\begin{equation} \label{eq:fnwre_mvp}
  \frac{\ntwo{Bx - \bar B x}}{\nfro{B} \ntwo{x}} \le \varepsilon.
\end{equation}
If $B$ is square and nonsingular (as it is in this application), then we can
obtain a simple expression for the maximum perturbation to $x$ that yields the
same error as $\bar B$ makes. Let $y = B x$. Let $\bar y \equiv \bar B x$; we
want to bound $\ntwo{\delta x}$ such that $\bar y = B (x + \delta
x)$. Rearranging, $\delta x = B^{-1} (\bar y - B x) = B^{-1} (\bar y - y)$. By
\eqref{eq:fnwre_mvp}, $\ntwo{\bar y - y} \le \varepsilon \nfro{B} \ntwo{x}$ and
so $\ntwo{\delta x} \le \varepsilon \nfro{B^{-1}} \nfro{B} \ntwo{x}$, or
\begin{equation} \label{eq:fnwre_mvp1}
  \frac{\ntwo{\delta x}}{\ntwo{x}} \le \varepsilon \kappa_{\text{F}}(B),
\end{equation}
where $\kappa_{\text{F}}(B)$ is the condition number of $B$ in the Frobenius
norm. If \eqref{eq:fnwre} instead has the form $\nfro{E} \le \varepsilon
\ntwo{B}$, then `2' replaces `F' in \eqref{eq:fnwre_mvp} and
\eqref{eq:fnwre_mvp1}.

One can also specify an element-wise error (EWE) of the approximation $\bar B$
to $B$. Let $\max_{i,j} |E_{i,j}| \equiv \nmax{E}$: this is the element-wise
max-norm. In this case, it is requested that $\bar B$ satisfy
\begin{equation} \label{eq:ewe}
  \nmax{E} \le \varepsilon N^{-1} \nrmone{B} \equiv \bar\varepsilon.
\end{equation}
If \eqref{eq:ewe} holds, then $\nrmone{Ex} \le \bar\varepsilon \nrmone{x}
\nrmone{e} = \bar\varepsilon \nrmone{x} N = \varepsilon \nrmone{B} \nrmone{x}$,
where $e$ is the vector of all ones, and so the relative error in an MVP is
bounded as
\begin{equation*}
  \frac{\nrmone{Bx - \bar B x}}{\nrmone{B} \nrmone{x}} \le \varepsilon,
\end{equation*}
which is the same as \eqref{eq:fnwre_mvp} except the 1-norm is used. Similarly,
one obtains the analogue of \eqref{eq:fnwre_mvp1}:
\begin{equation*}
  \frac{\nrmone{\delta x}}{\nrmone{x}} \le \varepsilon \kappa_1(B).
\end{equation*}

Observe that an element-wise norm (Frobenius, max) applied to blocks yields an
error bound in $p$-norms (respectively 2 and 1).

\section{Meeting the requested error tolerance} \label{sec:meeting}
Let a matrix be partitioned into blocks indexed by $i \in {\cal I}$; we often
write $\sum_i$ rather than $\sum_{i \in {\cal I}}$. We need to satisfy the bound
\eqref{eq:fnwre} on the whole matrix based on information in each block. The
bound may be achieved by at least three methods.

The first method---really, a class of methods---is to prescribe the rank of each
block based on analysis of problem-specific information: for example, an
integral equation and its discretization \citep{panel-cluster}; or in the
context of proving error bounds for far more sophisticated hierarchical
algorithms than we consider: for example, the {\it H$^2$}-matrix method of
\citet{borm-h2-variable}. Our interest is in applying the relatively simple
\hmatrix{} method to arbitrary kernels, and so we do not further consider these
approaches.

We call the second method the \emph{block-wise relative-error method}
(BREM). This method is standard in practice: for example, equation (3.69) and
surrounding text in \citet{bebendorf_book}, the software package AHMED by the
same author, numerical experiments in Chapter 4 of \citet{rjasanow-book}, and
relevant experiments in Section 4.8 of \citet{hlib} all use or describe this
method. Each LRA $\bar B_i$ to block $B_i$ is computed to satisfy
\begin{equation} \label{eq:bwrem}
  \nfro{E_i} \le \varepsilon \nfro{B_i}.
\end{equation}
Then
\begin{equation} \label{eq:bwrem-proof}
  \nfro{E}^2 = \sum_i \nfro{E_i}^2 \le \varepsilon^2 \sum_i \nfro{B_i}^2
  = \varepsilon^2 \nfro{B}^2,
\end{equation}
which implies \eqref{eq:fnwre}.

The third method supposes $\nfro{B}$ is known. We call it the \emph{matrix-wise
  relative-error method} (MREM). We believe that this method, though simple, is
new. Each LRA $\bar B_i$ to $B_i \in \mathbb{R}^{m_i \times n_i}$ is computed to
satisfy
\begin{equation} \label{eq:mwrem}
  \nfro{E_i} \le \varepsilon \frac{\sqrt{m_i n_i}}{N} \nfro{B}.
\end{equation}
As $N^2 = \sum_{i \in {\cal I}} m_i n_i$,
\begin{equation} \label{eq:mwrem-proof}
  \nfro{E}^2 = \sum_i \nfro{E_i}^2 \le \varepsilon^2 N^{-2} \nfro{B}^2 \sum_i
  m_i n_i = \varepsilon^2 \nfro{B}^2,
\end{equation}
which again implies \eqref{eq:fnwre}.

If \eqref{eq:bwrem} and \eqref{eq:mwrem} are equalities rather than
inequalities, then so are respectively \eqref{eq:bwrem-proof} and
\eqref{eq:mwrem-proof}.

The efficiencies of BREM and MREM are determined by the magnitudes $\nfro{B_i}$
and $N^{-1} \sqrt{m_i n_i} \nfro{B}$, respectively, as a function of $i \in
{\cal I}$. One magnitude cannot dominate the other for all $i \in {\cal I}$: as
equality is possible in \eqref{eq:bwrem-proof} and \eqref{eq:mwrem-proof}, if
for one block the first magnitude is greater than the other, then there is at
least one other block for which the opposite holds.

MREM requires that a block LRA satisfy (on rearranging \eqref{eq:mwrem})
\begin{equation} \label{eq:fnpe}
  (m_i n_i)^{-1} \nfro{E_i}^2 \le \varepsilon N^{-2} \nfro{B}^2.
\end{equation}
On each side of this inequality---omitting $\varepsilon$ on the right side---is
a quantity we call the \emph{Frobenius norm squared per element} (FNPE). If we
also square and divide \eqref{eq:bwrem} by $m_i n_i$, then our two magnitudes
from before become respectively the \emph{block} and \emph{matrix} FNPE $(m_i
n_i)^{-1} \nfro{B_i}^2$ and $N^{-2} \nfro{B}^2$. For a matrix $B$ arising from a
singular kernel, the block FNPE is smaller for far-off-diagonal blocks than for
near- and on-diagonal blocks. Hence relative to BREM, MREM requests less
accuracy for far-off-diagonal blocks and more for the others.

BREM and MREM require different termination criteria in the block LRA
algorithm. BREM requires the LRA algorithm to terminate based on an estimate of
the relative error with tolerance $\varepsilon$; MREM, the absolute error with
tolerance $\varepsilon N^{-1} \sqrt{m_i n_i} \nfro{B}$. Neither termination
criterion is more difficult to implement, or requires more work to evaluate
(within a small number of operations), than the other.

So far we have discussed the error bound in the Frobenius norm; now we discuss
the bound in the max-norm. If every block $i$ satisfies
\begin{equation} \label{eq:bewe}
  \nmax{E_i} \le \bar\varepsilon,
\end{equation}
then \eqref{eq:ewe} holds. We call this procedure MREMmax. Consider the
inequalities \eqref{eq:fnpe} and \eqref{eq:bewe}. In both, the right side is an
absolute tolerance that is the same for every block; and the left side describes
an element-wise quantity: in the first, the FNPE; in the second, the maximum
magnitude. Hence MREM and MREMmax behave similarly in how they map the matrix
tolerance to the block ones. It is not clear to us whether there are any
applications that use MREMmax rather than BREM. We shall shortly discuss why
MREM is preferable to MREMmax in practice.

\subsection{Estimating $\nfro{B}$} \label{sec:estnfro}
In some problems $\nfro{B}$ may be available, but we also need a means to
estimate this norm. We describe two methods.

The first is a stochastic method. Let $\{X_i\}$ be $n$ iid samples. The
estimator of the mean is, as usual, $\mu \equiv n^{-1} \sum_i X_i$. It is useful
to estimate confidence intervals on the estimator $\mu$ using a resampling
method. A straightforward approach is to compute the \emph{delete-1 jackknife
  variance} of the estimator, which in this simple case is
$n^{-1} \sum_i (X_i - \mu)^2$. We use the square root of this value and call it
the \emph{jackknife standard deviation} (JSD).

To estimate $\nfro{B}^2$, we could select a subset of entries. In practice we
choose to select a random subset of columns. This choice implies $X_i \equiv N
\sum_{j=1}^N B_{ji}^2$. We increase the number of columns until the JSD is less
than a requested tolerance.

A second method is to obtain an initial approximation $\tilde B$ to $B$ with a
very large tolerance $\tilde{\varepsilon}$. One can use either BREM or MREM; if
MREM, use the stochastic method to obtain the initial estimate of $\nfro{B}$. As
$\nfro{B - \tilde{B}} \le \tilde{\varepsilon} \nfro{B}$ and $\tilde{B} = B +
(\tilde{B} - B)$, $\nfro{\tilde{B}} \le (1 + \tilde{\varepsilon})
\nfro{B}$. Hence we can safely use $(1 + \tilde{\varepsilon})^{-1}
\nfro{\tilde{B}}$ as the estimate of $\nfro{B}$ in MREM to obtain the final
approximation $\bar B$. Computing $\nfro{\tilde{B}}$ requires work proportional
to $\nnz(\tilde{B})$.

\subsection{Recompression} \label{sec:recom}
Let $B \in \mathbb{R}^{m \times n}$ now be a block; in this subsection we
suppress the subscript $i$. Let $\bar B^1 = U^1 (V^1)^T$ be the rank-$r$ output
of an LRA algorithm. One can improve a suboptimal LRA by using a
\emph{recompression} algorithm \citep{bebendorf_book}. Such an algorithm
attempts to compress $\bar B^1$; the algorithm's efficiency results from the
outer-product structure of $\bar B^1$ and that $B$ is not accessed at all. For
the latter reason, the algorithm's error bound must be based on only $\bar
B^1$. An example of a recompression algorithm is to compute the singular value
decomposition (SVD) of $\bar B^1$---which is efficient because of the
outer-product structure of $\bar B^1$---and then to discard some subset of the
singular values and vectors. Let $\bar B^2$ be the output of a recompression
algorithm.

Let $\varepsilon$ be a block \emph{absolute} tolerance, as in MREM and
MREMmax. We must assure $\nfro{B - \bar B^2} \le \varepsilon$; what follows also
holds if `F' is replaced by `max'. Let $0 < \alpha < 1$. First, compute $\bar
B^1$ so that $\nfro{B - \bar B^1} \le \alpha \varepsilon$. Second, compute $\bar
B^2$ so that $\nfro{\bar B^1 - \bar B^2} \le (1 - \alpha) \varepsilon$; for
generality later, let $\beta \equiv 1 - \alpha$. As $B - \bar B^2 = (B - \bar
B^1) + (\bar B^1 - \bar B^2)$, $\nfro{B - \bar B^2} \le \nfro{B - \bar B^1} +
\nfro{\bar B^1 - \bar B^2} \le \alpha \varepsilon + (1 - \alpha) \varepsilon =
\varepsilon$, and so $\bar B^2$ satisfies the tolerance.

Now let $\varepsilon$ be a block \emph{relative} tolerance, as in BREM. We must
assure $\nfro{B - \bar B^2} \le \varepsilon \nfro{B}$. Again, let $0 < \alpha <
1$; we must determine $\beta > 0$. First, compute $\bar B^1$ so that $\nfro{B -
  \bar B^1} \le \alpha \varepsilon \nfro{B}$. Second, compute $\bar B^2$ so that
$\nfro{\bar B^1 - \bar B^2} \le \beta \varepsilon \nfro{\bar B^1}$. If in this
second inequality the right side involved $\nfro{B}$ rather than $\nfro{\bar
  B^1}$, then the calculation would be the same as in the case of absolute
error. But here we must bound $\nfro{\bar B^1}$. As $\bar B^1 = (\bar B^1 - B) +
B$, $\nfro{\bar B^1} \le \nfro{B - \bar B^1} + \nfro{B}$. From the first step,
$\nfro{B - \bar B^1} \le \alpha \varepsilon \nfro{B}$, and so $\nfro{\bar B^1}
\le (1 + \alpha \varepsilon) \nfro{B}$. Let $\beta = (1 - \alpha) / (1 + \alpha
\varepsilon)$. Then $\nfro{B - \bar B^2} \le \nfro{B - \bar B^1} + \nfro{\bar
  B^1 - \bar B^2} \le \alpha \varepsilon \nfro{B} + \beta \varepsilon \nfro{\bar
  B^1} \le \alpha \varepsilon \nfro{B} + (1 - \alpha) \varepsilon \nfro{B} =
\varepsilon \nfro{B}$, and so again $\bar B^2$ satisfies the tolerance.

In the case of SVD recompression, following Section 1.1.4 of
\citet{bebendorf_book}, let $Q_U R_U$ be the thin QR factorization of $U^1$ and
similarly for $V^1$. Let $W \Sigma Z^T$ be the SVD of $R_U R_V^T$. Then $Q_U W
\Sigma Z^T Q_V^T$ is the thin SVD of $\bar B^1$. This SVD is obtained with $O((m
+ n)r^2 + r^3)$ work. Choose the first $k \le r$ columns of $W$ and $V$ and
singular values $\sigma_i$ and set $U^2 = Q_U W_{\matcol \matcom 1 \matcol k}
\Sigma_{1 \matcol k \matcom 1 \matcol k}$ and $V^2 = Q_V Z_{\matcol \matcom 1
  \matcol k}$. This gives the LRA $\bar B^2 = U^2 (V^2)^T$; $\nnz(\bar B^2) \le
\nnz(\bar B^1)$, with strict inequality if $k < r$. When one uses MREM and so
the Frobenius norm, one chooses $k$ based on the inequality $\nfro{\bar B^2_k -
  \bar B^1} \le \beta \varepsilon$, where $\bar B^2_k$ uses the $k$ largest
singular values and associated vectors. Because $\nfro{\bar B^2_k - \bar B^1} =
(\sum_{i=k+1}^r \sigma_i^2)^{1/2}$, choosing $k$ requires just $O(r)$
operations. Using MREMmax entails more work because there is no simple
relationship between the singular values and the max-norm. One performs a binary
search on the ordered list of singular values. At each step, all the elements of
$\bar B^2_k$ must be computed and compared with $\bar B^1$ to find the maximum
deviation. Hence selecting $k$ requires $O(m n \log r)$ work; this work can
dominate the work to compute the SVD if $r$ is small relative to $m$ and
$n$. Because SVD recompression is a very effective way to improve an LRA, MREM
is preferable to MREMmax if one does not prefer the 1- to the 2-norm in the
associated error bounds.

\section{Numerical experiments}
We have several goals for our numerical experiments. First, of course, we must
show that MREM can yield greater compression than BREM on at least some problems
of interest. Second, we want to show that MREM is robust: ideally, it should
never yield \emph{less} compression than BREM. Third, we are interested in the
errors that MREM achieves: are they indeed only slightly better than requested?

We use \emph{adaptive cross approximation} (ACA) \citep{bebendorf-aca} as
implemented in AHMED to find block LRA. To implement the absolute block
tolerance required by MREM, we modified two lines of code. A good LRA algorithm
must compress a matrix well while requesting as few matrix entry evaluations as
possible. We find that ACA is quite efficient and robust on the problems in our
test set. It terminates when the relative error is within a factor of
approximately 10 either side of the tolerance; in particular, note carefully
that occasionally (but infrequently) the error is \emph{worse} than
requested. Of course ACA could terminate with a more precise error, but
increased termination precision would entail additional work. Instead, we
accommodate this behavior as follows. Let $\varepsilon$ be the block
tolerance. In the recompression tolerance we set $\alpha = 1/2$; if MREM is
used, $\beta = \alpha$; if BREM, $\beta$ is as described in Section
\ref{sec:recom}. We run ACA with a tolerance of $\alpha \varepsilon / 10$, where
the factor of $1/10$ compensates for ACA's sometimes imprecise termination
behavior, and then the SVD recompression algorithm with $\beta \varepsilon$. LRA
algorithms based on interpolation ({\it e.g.}, \citet{hybrid-aca}) mitigate the
problem of terminating with less accuracy than requested while still remaining
efficient. In any case, the choice of LRA algorithm is independent of the choice
of either BREM or MREM.

We estimate $\nfro{B}$ using the stochastic method described in
\ref{sec:estnfro}. The estimator is terminated when the JSD is no greater than
$1/50$. To be conservative, we subtract two JSD from the estimate of $\nfro{B}$,
thereby effectively decreasing the tolerance slightly.

\subsection{Particle distributions}
We perform numerical experiments using a distribution of point-like particles
that interact according to the kernel $K(r) = 1/r^p$ if $r > 0$, $0$ otherwise,
where $r$ is the distance between two points and $p$ is the order of the
singularity. We also consider the kernel $K(r) = \ln r$.

We consider singularity orders $p = 1$, $2$, $3$ and three geometries: uniform
distribution of points in the cube $[-1,1]^3$, on its surface, and on its edges;
where indicated in figures, these are denoted respectively `cube', `surf', and
`edge'. The requested error tolerance $\varepsilon$ is usually $10^{-5}$. We
measure several quantities. The one that most clearly demonstrates the
performance of MREM relative to BREM is the \emph{improvement factor} (IF),
$\nnz(\bar B^{\text{B}})/\nnz(\bar B^{\text{M}})$, where the superscripts B and
M denote BREM and MREM. A value greater than one indicates improvement.

\begin{figure}
\centering
\includegraphics[width=5in]{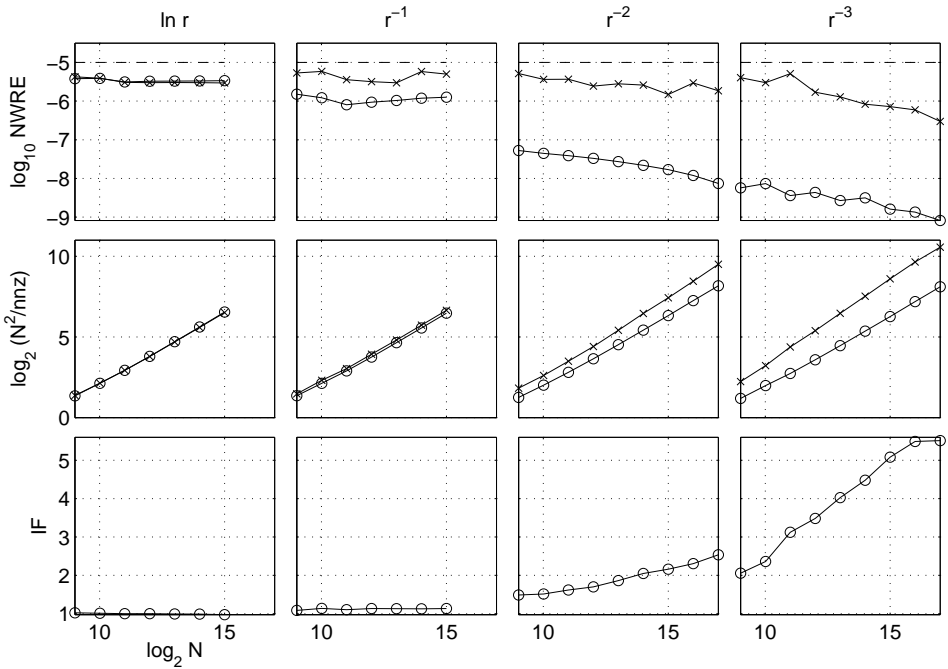}
\caption{Particles are distributed uniformly on the edges of a cube.}
\label{fig:Ged}
\end{figure}

\begin{figure}
\centering
\includegraphics[width=5in]{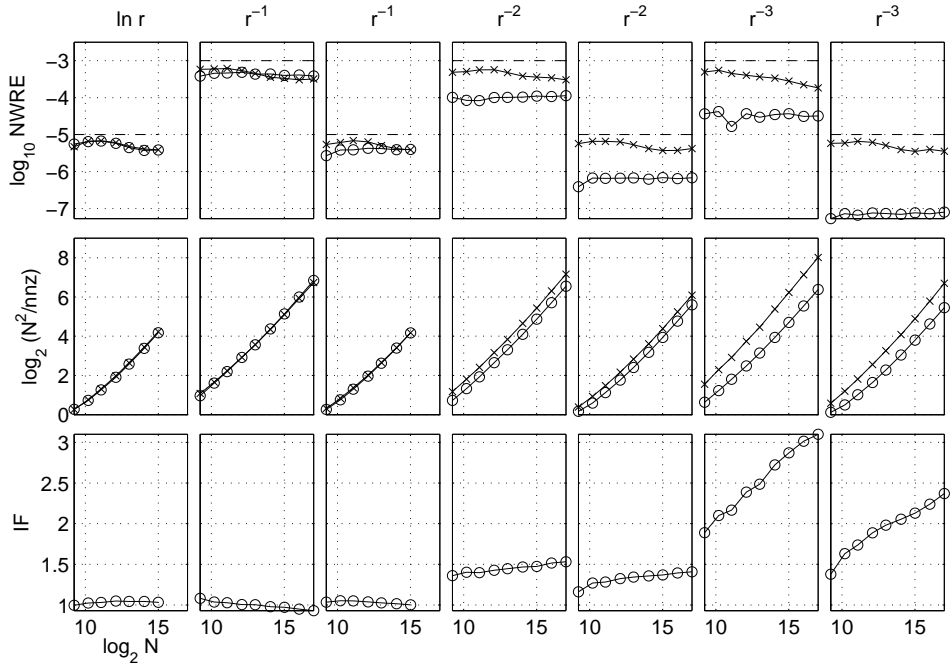}
\caption{Particles are distributed uniformly on the surface of a cube.}
\label{fig:Gsu}
\end{figure}

\begin{figure}
\centering
\includegraphics[width=5in]{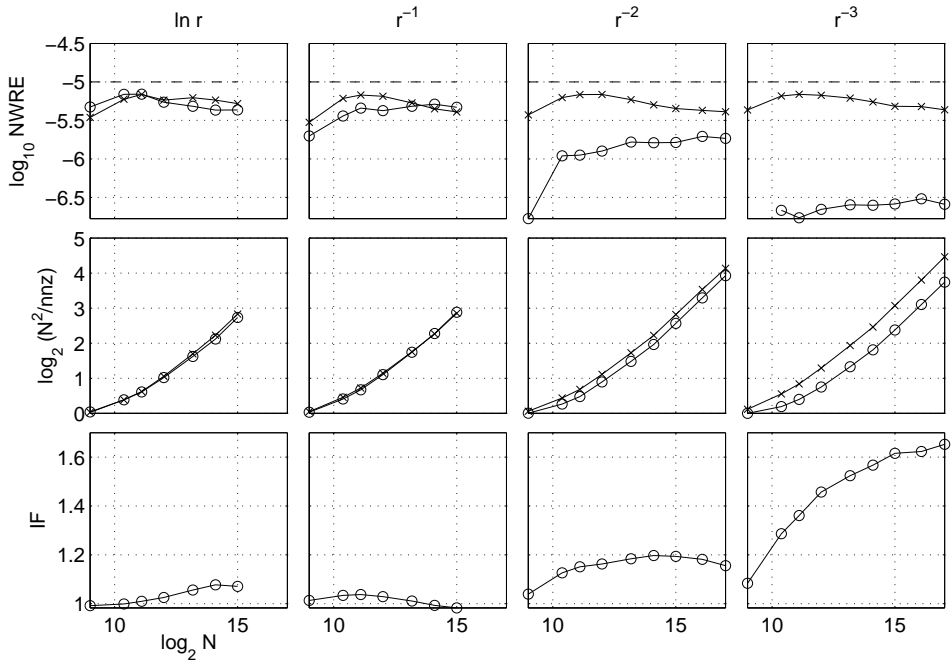}
\caption{Particles are distributed uniformly inside a cube.}
\label{fig:Gcu}
\end{figure}

Figures \ref{fig:Ged}, \ref{fig:Gsu}, and \ref{fig:Gcu} show results for the
three geometries. Each column corresponds to a kernel, which is stated at the
top of the column, and tolerance $\varepsilon$, indicated by the dashed
horizontal line in the top row. In the top two rows, curves for BREM have
circles; for MREM, {\it x}'s. All rows plot quantities against problem size
$N$. The top row shows the Frobenius-norm-wise relative error achieved; these
are estimates for problems having size $N > 2^{15}$. The middle row shows the
compression factor. The bottom row shows the improvement factor of MREM. Most
experiments were carried out to approximately $2^{17}$ particles, though in a
few cases the runs were terminated early to save computation time or because of
memory size.

Trends accord with our observations in Section \ref{sec:meeting}. The IF
increases as the dominance of the (near) diagonal of $B$ increases. In the
context of our test problems, the IF increases with increasing order $p$ and
sparser geometry. The IF rarely is below 1, and then by only an extremely small
amount and only on problems whose singularity is of order $1$. For both BREM and
MREM, achieved errors are always below the requested tolerance; and MREM's
achieved errors are almost always within a factor of $10$ of it.

\begin{figure}
\centering
\includegraphics[width=5in]{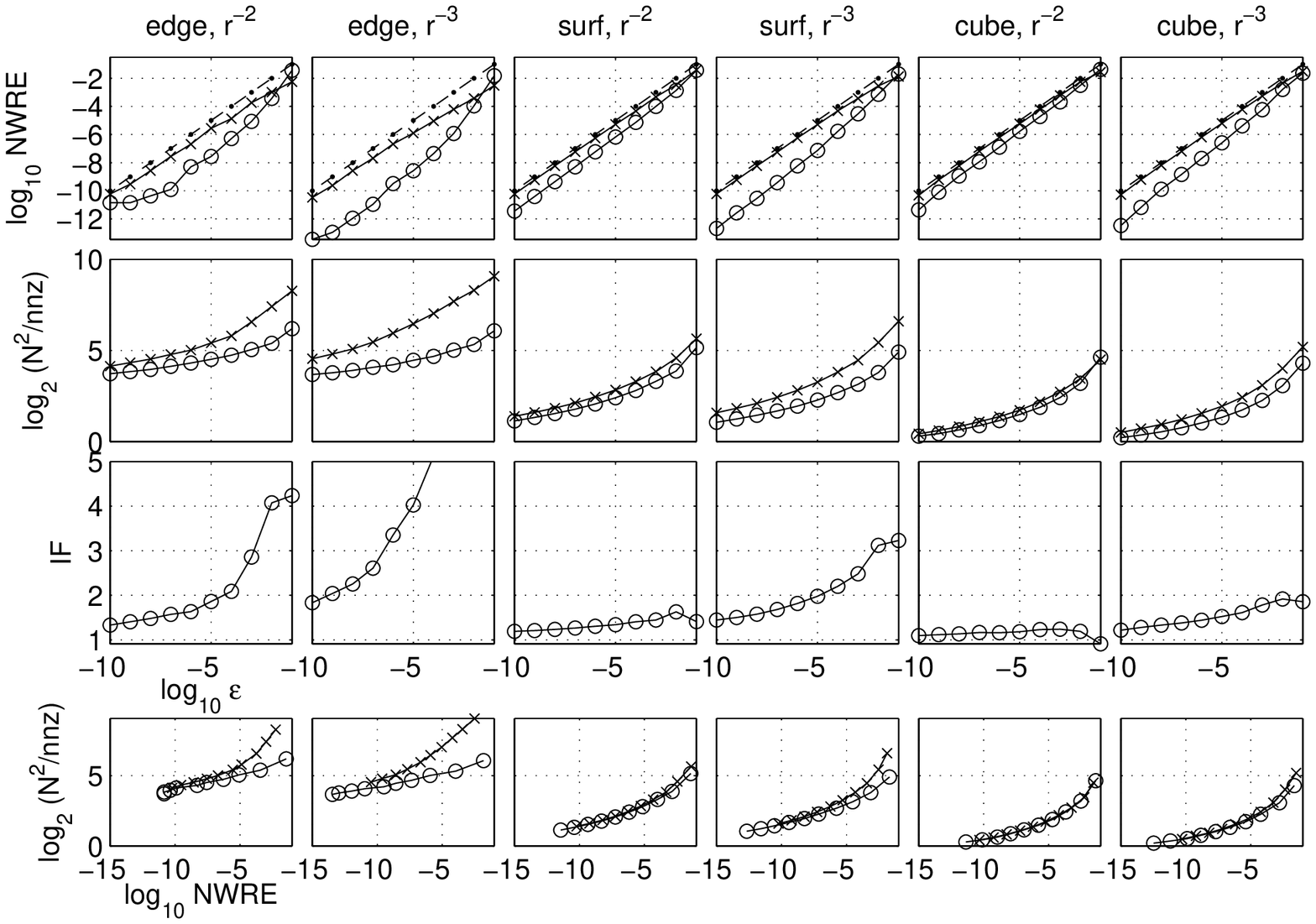}
\caption{The requested tolerance $\varepsilon$ is varied. The number of
  particles is approximately $N = 2^{13}$.}
\label{fig:tol}
\end{figure}

Figure \ref{fig:tol} shows the behavior of the two methods for $p = 2,3$ and the
three geometries when the requested tolerance is varied. The first three rows
are as in the earlier figures except that the abscissa is the tolerance rather
than problem size. For easy reference, in the first row requested tolerance is
also indicated by the line having dots. The fourth row plots compression against
\emph{achieved} error. On these test problems, MREM achieves greater compression
than BREM even when the achieved errors are the same.

\subsection{A problem in geophysics}
The problem that motivated this work is one in geophysics. We are modeling the
evolution of slip on a fault in an elastic half space. One part of the
computation is a matrix-vector product relating slip on the fault to stress. The
fault is discretized by rectangular patches, and stress at the center of all
patches due to unit slip on one patch is computed using the rectangular
dislocation code of \citet{okada92}. The matrix is approximated by an
\hmatrix{}. We performed a test on a fault discretized by $156 \times 402$
rectangles with a tolerance of $10^{-5}$. For respectively BREM and MREM, the
actual errors are $1.28 \times 10^{-8}$ and $1.71 \times 10^{-6}$; the
compression factors are $16.89$ and $75.74$, for an improvement factor of
$4.48$.

\section{Conclusions}
For many problems, MREM produces a more efficient approximation than BREM for a
requested error tolerance by producing an approximation $\bar B$ that is little
more accurate than is requested. MREM rarely produces a less efficient
approximation, and then only by a small amount. Improvement factors are often
between $1.5$ and $5$. The improvement factor increases with the dominance of
the diagonal. For problems in which the order of the singularity is $1$, MREM's
results are almost exactly the same as BREM's. For higher-order problems, MREM
is consistently better by a substantial and practically useful amount. MREM is
also better than BREM for high-order singularities even when \emph{achieved},
rather than requested, error is considered.

Using MREM rather than BREM requires the extra work to estimate $\nfro{B}$ if it
is not known. In our experiments, $\nfro{B}$ is estimated in a time that is
approximately $1\%$ of the total construction time. Because the work to find a
block LRA increases with the block tolerance, MREM is faster than BREM on any
problem in which the resulting compression is greater by at least a small
amount.

The similarity in block tolerances suggests that the compression factors
resulting from MREM and MREMmax can be expected to scale with problem size and
tolerance similarly. However, SVD recompression is far more efficient when using
MREM. We conclude that one should consider using MREM when a problem's
singularity has order greater than $1$. \\

\vspace{-.2cm} {\bf Acknowledgements.} I would like to thank Prof.~Eric Darve,
Dept.~of Mech.~Eng., Stanford University, for reading a draft of this paper and
a helpful discussion; and Prof.~Paul Segall, Dept.~of Geophysics, Stanford
University, for supporting this research.  \bibliographystyle{siam}
\bibliography{../amb_cdfm}{}

\end{document}